\documentclass[12pt]{article}
\begin{document}

\title{Hirzebruch genera of manifolds equipped with a Hamiltonian circle 
action}
\author{K. E. Feldman\footnote{The author was supported by 
Seggie Brown fellowship.}}
\date{}
\maketitle

In the paper we describe obstructions for the existence of symplectic and 
Hamiltonian symplectic circle actions on closed compact manifolds in terms of 
Hirzebruch genera and relations between differential and homotopic invariants 
of such manifolds. All manifolds considered are assumed to be closed and 
compact.

\

{\bf Theorem 1.} The Todd genus of a manifold equipped with a symplectic 
circle action with isolated fixed points is either equal to zero and then
the action is non-Hamiltonian, or equal to one and then the action is 
Hamiltonian. Any symplectic circle action on a manifold with the positive 
Todd genus is Hamiltonian.

{\bf Proof.} 
A symplectic circle action is Hamiltonian if and only if there is such
a connected 
component of the fixed point set that its  weights of the $S^1$-representation
 in the normal bundle are positive. This connected component is unique 
and corresponds to the global minimum of the moment map (see details 
in~\cite{Audin}). Let us apply the fixed point formula for Hirzebruch genus 
$\chi_y$ corresponding to the power series 
$\frac{x\left(1+ye^{-x(1+y)}\right)}{1-e^{-x(1+y)}}$.
Let $M_s$ be connected components of the fixed point set, $d_s$ be the number 
of negative weights of the $S^1$-representation in the normal bundle 
of $M_s$, then(see~\cite{HBJ})
\begin{equation}
\label{1}
\chi_y(M)=\sum_s(-y)^{d_s}\chi_y(M_s),
\end{equation}
here the sum is taken over all connected components of the fixed point set.
Since $\chi_y(pt)=1$, the constant term of  $\chi_y$ of a manifold 
equipped with a Hamiltonian circle action with isolated fixed points is equal 
to one. Vice versa, if the constant term of $\chi_y(M)$ is not equal to zero, 
then any symplectic circle action on  $M$ is Hamiltonian. The statement of 
the theorem follows from $Td(M)=\chi_y(M)|_{y=0}$.

\

The result of Theorem 1 is a particular case of the restrictions imposed 
on manifolds with a circle action by the Conner--Floyd equations.
Let us fix a multiplicative genus  $h:\Omega^U_*\to Q$.
Let $g_h(u)=u+\sum\frac{h([CP(n)])}{n+1} u^{n+1}$. 
For an almost complex manifold $M^{2n}$, equipped with a quasi-complex circle 
action with isolated fixed points, the sum over all the fixed points of the 
products
$\prod^n_{k=1}\left(g^{-1}_h (j_k g(u))\right)^{-1}$($j_k$, $k=1,...,n$, 
are weights of the $S^1$-representation in a neighborhood of the 
fixed point considered)
contains no negative powers of $u$ and its constant term is equal to 
$h(M^{2n})$ (see~\cite{Novikov,Krichever}).

\

{\bf Corollary 1\cite{TolmanWeitsman}.} A semi-free symplectic circle action 
with isolated fixed points is Hamiltonian if and only if the fixed point set 
is non-empty.

{\bf Proof.} Let us write the Conner-Floyd equations for the Todd genus in 
this case. Because $Td([CP(m)])=1$ for every $m$, we obtain 
$g_{Td}(u)=\sum u^{m+1}/(m+1)=ln(1-u)$, $g^{-1}_{Td}(jg_{Td}(u))=1-(1-u)^j$.
All the weights of the $S^1$-representation are equal to $\pm 1$. 
Summing up over all the fixed points we obtain that the function
$f(u)=\sum^{2n}_{k=1} F_k\frac{(u-1)^k}{u^n}$,
where $F_k$ is the number of fixed points with exactly $k$ negative weights,
is a constant, say $\lambda$, which is equal to the Todd 
genus of the manifold. Solving the corresponding system of linear with 
respect to $F_k$ equations, we find that
$(F_0,F_1,...,F_n)=\lambda\left(\left(n\atop 0\right), \left(n\atop 1\right),...,\left(n\atop n\right)\right)$.
Applying Theorem 1 we deduce the statement of the corollary.

{\bf Corollary 2\cite{McDuff}.} A symplectic circle action on a 
$4$-dimensional manifold is Hamiltonian if and only if the fixed point set is 
non-empty.

{\bf Proof.} Let us show that the existence of a fixed point implies the 
existence of a fixed point whose weights of the $S^1$-representation 
are of the same sign. Otherwise there is no connected component of 
codimension 2 (such component has only one weight). Thus, all fixed points 
are isolated and have exactly one negative and one positive weights. 
The contradiction follows from:

{\bf Lemma 1.} Every quasi-complex circle action with isolated fixed points on 
$M^{2n}$ has fixed points with different parity of the number of negative 
weights. 

{\bf Proof.} Otherwise the coefficients of $1/u^n$ in the Conner-Floyd 
equations for the Todd genus of the manifold have the same sign.

\

 Let us consider in detail the case of a Hamiltonian circle action. Denote
 the Poincare polynomial of $M$ by $P_M(t)$.

{\bf Theorem 2.} If a symplectic manifold $M$ admits a Hamiltonian circle 
action with isolated fixed points, then $\chi_y(M)=P_M(\sqrt{-y})$.

{\bf Proof.} Critical points of the moment map coincide with the fixed points 
of the circle action. The index of every critical point is twice the number of 
negative weights of the $S^1$-representation. Thus, the statement of  
the theorem follows from~(\ref{1}).

{\bf Corollary 3\cite{JonesRawnsley}.} The signature and the Betti numbers of 
a manifold $M^{4n}$ equipped with a symplectic circle action with isolated 
fixed points are related through the formula:
$\sigma(M^{4n})=\sum^n_{k=0} b_{4k} -\sum^{n-1}_{k=0} b_{4k+2}$.

{\bf Proof.} The statement follows from Theorem 2 after substitution
$y=1$.

{\bf Corollary 4.} Homotopy equivalent manifolds equipped with Hamiltonian 
circle actions with isolated fixed points have the same $\chi_y$-genus. 

{\bf Remark 1.} Similar to~\cite{Farber} for a non-Hamiltonian 
symplectic circle action with isolated fixed points 
$\chi_y(M^{2n})=\sum^n_{k=0} b_{2k}(-y)^k$, here $b_k$ are Novikov 
numbers~\cite{Novikov1} corresponding to the generalized moment map.
It is important to construct an example of such an action with a
non-empty fixed point set.

\

University of Edinburgh

feldman@maths.ed.ac.uk

\end{document}